\documentclass[12pt]{article}
\usepackage[margin=2cm]{geometry}
\usepackage{amsfonts}
\usepackage{amssymb}
\usepackage{amsthm}
\usepackage{amsmath}
\usepackage{latexsym}
\usepackage{color}
\usepackage{mathrsfs}

\begin{document}
\allowdisplaybreaks[4]
\newtheorem{lemma}{Lemma}
\newtheorem{pron}{Proposition}
\newtheorem{thm}{Theorem}
\newtheorem{Corol}{Corollary}
\newtheorem{exam}{Example}
\newtheorem{defin}{Definition}
\newtheorem{remark}{Remark}
\newtheorem{property}{Property}
\newcommand{\la}{\frac{1}{\lambda}}
\newcommand{\sectemul}{\arabic{section}}
\renewcommand{\theequation}{\sectemul.\arabic{equation}}
\renewcommand{\thepron}{\sectemul.\arabic{pron}}
\renewcommand{\thelemma}{\sectemul.\arabic{lemma}}
\renewcommand{\thethm}{\sectemul.\arabic{thm}}
\renewcommand{\theCorol}{\sectemul.\arabic{Corol}}
\renewcommand{\theexam}{\sectemul.\arabic{exam}}
\renewcommand{\thedefin}{\sectemul.\arabic{defin}}
\renewcommand{\theremark}{\sectemul.\arabic{remark}}
\renewcommand{\theproperty}{\sectemul.\arabic{property}}
\def\REF#1{\par\hangindent\parindent\indent\llap{#1\enspace}\ignorespaces}

\title{\large\bf On the random-time and finite-time ruin probability for widely dependent claim sizes and inter-arrival times}
\author{\small Yang Chen$^1$\thanks{Research supported by Humanities and Social Sciences Foundation of the Ministry of Education of China (No. 18YJC910004).},~~Zhaolei Cui$^2$~~and~~Yuebao Wang$^3$
\thanks{Corresponding author.
Telephone: +86 512 67422726. Fax: +86 512 65112637. E-mail:
ybwang@suda.edu.cn}
\\
{\footnotesize\it 1 School of Mathematical Sciences, Suzhou University of Science and Technology, Suzhou, 215009, China}\\
{\footnotesize\it 2 School of Mathematics and Statistics, Suzhou University of Technology, Suzhou 215000, China}\\
{\footnotesize\it 3 School of Mathematical Sciences, Soochow University, Suzhou 215006, China}\\
}
\date{}

\maketitle {\noindent\small {\bf Abstract }}\\

Using the results of precise large deviation and renewal theory for widely dependent random variables,
this paper obtains the asymptotic estimation of the random-time ruin probability
and the uniform asymptotic estimation of finite-time ruin probability for a nonstandard renewal risk model,
in which both claim sizes and the inter-arrival times of claim sizes are widely dependent.
\\



\section{\normalsize\bf Introduction}
\setcounter{equation}{0}\setcounter{thm}{0}\setcounter{lemma}{0}\setcounter{remark}{0}\setcounter{pron}{0}\setcounter{Corol}{0}
\setcounter{Corol}{0}

In a project of insurance business, we often use the finite-time ruin probability to estimate
the possibility of failure of the project within a given operating time.
However, under the influence of various uncertain factors, the operating time of the project also is uncertain or unknowable. In essence, the operating time frequently behaves as a random variable.
Therefore, when assessing the risk associated with an insurance business, it is imperative to examine both the finite-time ruin probability and the random-time ruin probability.
In addition, in the practice of insurance business, the participants are often not independent of each other,
that is, the actual risk model is often non-standard.
Therefore, we also need to consider some suitable dependent structures of random variables to describe the relationship between these participants, and analyze whether these dependent structures have an impact on the ruin probability.

In this paper, for a nonstandard renewal risk model with the widely dependent claim sizes and the widely dependent inter-arrival times of claim sizes, see Subsections 1.2 and 1.3 below, we give the asymptotic estimation of the corresponding random-time ruin probability
and the uniform asymptotic estimation of finite-time ruin probability for some range of operating time, see Theorems \ref{thm104}-\ref{thm106} below. The consistency in the asymptotic estimation of the finite-time ruin probability demonstrates that, within a specific range, the operating time of an insurance business is independent of its initial capital in a sense.
Consequently, this estimation holds greater theoretical significance (including technical difficulty)
and practical value compared to the corresponding non-uniform estimation.

In the existing literature, most studies first establish the uniform asymptotic estimation of finite-time ruin probability,
then directly derive corresponding asymptotic results of random-time ruin probability according to the former.
However, these established approaches typically impose the restrictive assumption that claim sizes must be mutually independent,
as seen in Theorems \ref{thm101}-\ref{thm103} below.
Departing from conventional methods, our study adopts a novel approach:
we first derive the asymptotic estimation for random-time ruin probabilities by applying results
on precise large deviations and renewal theory for widely dependent random variables.
In this way, we can get corresponding results for finite-time ruin probability directly,
then we study its uniform estimation of  over an operating time interval.
Both the proof line and methodology differ significantly from those found in prior works.

Moreover, it is well established that such results fundamentally depend
on the underlying distributional properties of the research objects and the relationship between them.
Accordingly, before presenting the main results of this paper,
we first introduce the relevant distribution classes,
the wide dependence structure of random variables, and a related nonstandard risk model
to provide the necessary theoretical foundation.

\subsection{\normalsize\bf Some distribution classes}

Throughout this paper, the following notations and conventions are adopted unless otherwise specified.

All limit relations refer to $x\to\infty$ and all distributions are supported on $(-\infty,\infty)$
(including $[0,\infty)$ or $(0,\infty)$).
Let $s_1(\cdot)$ and $s_2(\cdot)$ be two positive functions on $[0,\infty)$.
For $1\le i\neq j\le2$, we set
$$s_{i,j}=\limsup s_i(x)s_j^{-1}(x).$$
Then $s_i(x)=O\big(s_j(x)\big)$ means $s_{i,j}<\infty$, $s_1(x)\asymp s_2(x)$ means $\max\{s_{1,2},s_{2,1}\}<\infty$,
$s_i(x)\lesssim s_j(x)$ means $s_{i,j}\le1$, $s_1(x)\sim s_2(x)$ means $s_{1,2}=s_{2,1}=1$,
and $s_i(x)=o\big(s_j(x)\big)$ means \textcolor[rgb]{1.00,0.00,0.00}{$s_{i,j}=0$}.

Let $V$ be a distribution. We denote the tail of $V$ by $\overline{V}=1-V$,
and the $n$-fold convolution of $V$ with itself by $V^{*n}$ for all nonnegative integers $n$,
where $V^{*1}=V$ and $V^{*0}$ is the distribution of random variables degenerate to $0$.

We say that the distribution class
\begin{eqnarray*}
{\mathcal L}=\big\{V:\overline{V}(x-t)\sim\overline{V}(x)\ \text{for each}\ t\in(-\infty,\infty)\big\}
\end{eqnarray*}
is long-tailed, and the distribution class
\begin{eqnarray*}
{\mathcal S}=\big\{V:V\in\mathcal{L}\ \ \text{and}\ \ \overline{V^{*2}}(x)\sim2\overline{V}(x)\big\}
\end{eqnarray*}
is subexponential, which was introduced by Chistyakov \cite{C1964}.
Further, Lemma 2 of Chistyakov \cite{C1964} notes that the requirement $V\in\mathcal{L}$ is unnecessary if $V$ is supported on $[0,\infty)$ or $(0,\infty)$.
For a detailed discussion of classes $\mathcal{S}$ and $\mathcal{L}$, see Embrechts et al. \cite{EKM1997},
Resnick \cite{R2007}, Borovkov and Borovkov \cite{BB2008}, Foss et al. \cite{FKZ2013}, Wang \cite{W2022},
Leipus et al. \cite{LSK2023}, Wang et al. \cite{WCX2025}, and so on.

If $V\in\mathcal{S}$, then the following set of positive functions is not empty:
\begin{eqnarray}\label{101}
\mathcal{H}_V=\Big\{h(\cdot)\ \text{on}\ [0,\infty):h(x)\uparrow\infty,h(x)x^{-1}\downarrow0
\ \text{and}\ \int_{h(x)}^{x-h(x)}\overline{V}(x-y)V(dy)=o\big(\overline{V}(x)\big)\Big\}.
\end{eqnarray}

An important subclass of $\mathcal{S}$ introduced by Kl\"{u}ppelberg \cite{K1988} is
\begin{eqnarray*}
{\mathcal S^*}=\Big\{V:0<\int_0^\infty yV(dy)<\infty\ \ \text{and}\ \
\int_0^x\overline{V}(x-y)\overline{V}(y)dy\sim2\int_0^\infty yV(dy)\overline{V}(x)\Big\}.
\end{eqnarray*}
In order to introduce some other distribution classes, for each $y\in(0,\infty)$, we write
$$\overline{V_{*}}(y)=\liminf{\overline V}(xy){\overline V}^{-1}(x),
~~\overline{V^{*}}(y)=\limsup{\overline V}(xy){\overline V}^{-1}(x)\ \ \text{and}\ \ L_V=\lim\limits_{y\downarrow1}\overline{V_*}(y).$$
Then we say that the distribution classes
\begin{eqnarray*}
\mathcal{D}=\big\{V:{\overline{V_*}}(y)>0\ \ \text{for each}\ y\in(1,\infty)\big\},\ \ {\cal C}=\{V:L_V=1=\lim\limits_{y\uparrow1}\overline{V^*}(y)\},
\end{eqnarray*}
\begin{eqnarray*}
\mathcal{ERV}(\alpha,\beta)=\{V:y^{-\beta}\leq\overline{V_*}(y)\leq\overline{V^*}(y)\leq y^{-\alpha}
\ \ \text{for each}\ y>1\}
\end{eqnarray*}
for each pair $0\le\alpha\le\beta<\infty$ and
$$\mathcal{ERV}=\bigcup_{0\le\alpha\le\beta<\infty}\mathcal{ERV}(\alpha,\beta)$$
are dominated varying tailed, consistently varying tailed,
extended regularly varying tailed with indexes $0\le\alpha\le\beta<\infty$ and extended regularly varying tailed, respectively.
Particularly, if $\alpha=\beta$, then $\mathcal{ERV}(\alpha,\beta)$ reduces to the regularly varying tailed distribution class, denoted
by $\mathcal{R}_{\alpha}$.

Some properties are summarized into the following two propositions.

\begin{pron}\label{pron101}
The following inclusion relations are proper: for each pair $0<\alpha\le\beta<\infty$,
$$\mathcal{R}_{\alpha}\subset \mathcal{ERV}(\alpha,\beta)\subset\mathcal{C}\subset\mathcal{L}\cap\mathcal{D}\subset\mathcal{S}\subset\mathcal{L}.$$
\end{pron}

Further, we denote the moment index of $V$ by
$$I_V=\sup\Big\{s:\int_0^\infty y^sV(dy)<\infty\Big\},$$
and the upper Matuszewska index and lower Matuszewska index of distribution $V$ by
\begin{eqnarray*}
J_V^+=-\lim\limits_{y\to\infty}\ln \overline{V_*}(y)\ln^{-1} y\ \ \ \ \text{and}\ \ \ \
J_V^-=-\lim\limits_{y\to\infty}\ln \overline{V^*}(y)\ln^{-1} y.
\end{eqnarray*}

\begin{pron}\label{pron102} Let $V$ be a distribution in class $\mathcal{D}$.

$(1)$ $0\le J_V^-\le I_V\le J_V^+<\infty$.

$(2)$ For each $p<J_V^-$, $\overline{V}(x)=o(x^{-p})$.

$(3)$ For each $p>J_V^+$, $x^{-p}=o\big(\overline{V}(x)\big)$.
\end{pron}

The above concepts and properties can be found in some references,
such as Feller \cite{F1971}, Bingham et al. \cite{BGT1987}, Cline and Samorodnitsky \cite{CS1994}, Embrechts et al. $\cite{EKM1997}$,
Kl\"{u}ppelberg and Mikosch $\cite{KM1997}$ and Tang and Tsitsiashvili $\cite{TT2003}$.

\subsection{\normalsize\bf Wide dependent structure}

Wang et al. $\cite{WWG2013}$ introduced the concept of wide dependence structure of random variables.
By definition, random variables $\xi_i,i\ge 1$ with common distribution $V$ are said to be widely upper orthant dependent (WUOD),
if for each $n\ge1$, there exists some positive number $g_{U,V}(n)$ such that,
\begin{eqnarray}\label{102}
P\Big(\bigcap^{n}_{i=1}\{\xi_i>x_i\}\Big)\leq g_{U,V}(n)\prod_{i=1}^nP(\xi_i> x_i),\ \ \ \ \ \ x_i\in(-\infty,\infty),\ \ \ 1\le i\le n;
\end{eqnarray}
they are said to be widely lower orthant dependent (WLOD),
if for each $n\ge1$, there exists some  positive number $g_{L,V}(n)$ such that,
\begin{eqnarray}\label{103}
P\Big(\bigcap^{n}_{i=1}\{\xi_i \leq x_i\}\Big)\leq g_{L,V}(n)\prod_{i=1}^n P(\xi_i\leq x_i),
\ \ \ \ \ x_i\in(-\infty,\infty),\ \ \ 1\le i\le n;
\end{eqnarray}
and they are said to be widely orthant dependent (WOD) if they are both WUOD and WLOD.

WUOD, WLOD and WOD structures can be called  widely dependent (WD) as a joint name.
And $g_{U,V}(n),~g_{L,V}(n),\ n\geq 1$, are called dominating coefficients.
Clearly,
$$g_{U,V}(n+1)\ge g_{U,V}(n)\ge g_{U,V}(1)=1\ \ \ \text{and}\ \ \ g_{L,V}(n+1)\ge g_{L,V}(n)\ge g_{L,V}(1)=1,\ \ n\ge2.$$

Some basic properties of WD random variables are as follows, see Proposition 1.1 of Wang et al. $\cite{WWG2013}$.

\begin{pron}\label{Proposition103}
$(1)$ Let $\xi_i,i\geq1$ be WLOD (or WUOD) with common distribution $V$ on $(-\infty,\infty)$. If $f_i(\cdot),i\geq1$ are nondecreasing,
then $f_i(\xi_i),i\geq1$ are still WLOD (or WUOD);
if $f_i(\cdot),i\geq1$ are nonincreasing, then $f_i(\xi_i),i\geq1$ are WUOD (or WLOD).

$(2)$ If $\xi_i,i\geq1$ are nonnegative and WUOD with common distribution $V$ and finite mean,
then for each $n\geq1$,
$$E\prod^{n}_{i=1}\xi_{i} \leq g_{U,V}(n)\prod^{n}_{i=1}E\xi_{i}.$$
Therefore, if $\xi_i,i\geq1$ are WUOD with common distribution $V$ on $(-\infty,\infty)$ satisfying $Ee^{s {\xi_i}}<\infty$
for some $s>0$, then for each $n\geq1$,
$$Ee^{s \sum^{n}_{i=1}\xi_{i}} \leq g_{U,V}(n)\prod^{n}_{i=1}Ee^{s {\xi_i}}.$$
\end{pron}

Furthermore, Wang et al. \cite{WWG2013} demonstrated through examples that the WD structure encompasses various dependence types, including common negatively dependent random variables, certain positively dependent random variables,
and other dependence random variables.

When $g_{U,V}(n)=g_{L,V}(n)=M$ for all $n\ge1$ and some $M>0$,
the inequalities (\ref{101}) and (\ref{102}) describe extended negatively upper
and lower orthant dependent (ENUOD and ENLOD) random variables, respectively.
$\xi_i,i\ge 1$ are said to be extended negatively orthant dependent (ENOD) if they are both ENUOD and ENLOD.
ENOD, ENUOD, ENLOD random variables are collectively called END r.v.s, see Liu $\cite{L2009}$.

Further, if $M=1$, then we have the corresponding notions of NUOD, NLOD, NOD and ND random variables,
see, for example, Ebrahimi and Ghosh $\cite{EG1981}$ and Block et al. $\cite{BSS1982}$.

In particular, we say that random variables $\xi_i,i\ge1$ are negatively associated (NA),
if for any disjoint nonempty subsets $A$ and
$B$ of $\{1,\cdots,m\},\ m\geq2$ and any coordinate-wise
nondecreasing and positive functions $l_1(\cdot)$ and $l_2(\cdot)$, the inequality
\begin{eqnarray*}
Cov(l_1(\xi_i:i\in A),l_2(\xi_j:j\in B))\leq0
\end{eqnarray*}
holds whenever the moment involved exists.
For further details, please refer to Joag-Dev and Proschan \cite{JP1983} among others.

In real-world complex systems, random variables frequently exhibit dependence structures that go beyond strict independence or even negative dependence. For more information, please refer to the introduction before Proposition \ref{Proposition104} below.
Consequently, the wide dependence framework offers a more comprehensive and realistic modeling approach compared to models restricted to independence or negative dependence assumptions.



\subsection{\normalsize\bf Nonstandard renewal risk models}

Here, we introduce some concepts and properties of renewal risk model for the insurance business.

In this model, let claim sizes $Y_i,\ i\geq1$ be random variables
with common distribution $G$ on $[0,\infty)$ and mean $0<\mu_G<\infty$;
the inter-arrival times $Z_i,\ i\geq1$ of the claim sizes also be random variables
with common distribution $H$ on $[0,\infty)$ and mean $0<\mu_H<\infty$.

The times of successive claims, $S^H_n=\sum_{i=1}^nZ_i,\ n\geq1$, constitute a counting process
$$\Big(N(t)=\sum_{n=1}^\infty\textbf{1}_{\{S^H_n\leq\  t\}}:t\in[0,\infty)\Big)$$
with the corresponding mean function $\lambda(\cdot)=EN(\cdot)$.
For convenience, we might as well assume that
$$\inf\{t:0<\lambda(t)<\infty\}=0.$$
In addition, let $c$ be a positive constant interest rate and $x$ be a nonnegative initial capital of the insurance business.
In order to ensure normal operation of the insurance business, a safe load condition is usually assumed that
$$\mu_G<c\mu_H.$$

The aforementioned model is referred to as the standard renewal risk model,
if the random vectors $(Y_i,Z_i),i\geq1$ are independent and identically distributed,
and $\{Y_i,\ i\geq1\}$ are independent of $\{Z_i,\ i\geq1\}$.

Under these assumptions, we define $\big(N(t):t\in[0,\infty)\big)$ as the standard renewal counting process and $\lambda(\cdot)=EN(\cdot)$ as the standard renewal function.

Alternatively, if at least one of these independent assumptions is not true, the model is considered nonstandard.
In this case, $\big(N(t):t\in[0,\infty)\big)$ and $\lambda(\cdot)$ are respectively
called quasi-renewal counting process and quasi-renewal function.
They are collectively called the renewal risk model, the renewal counting process and the renewal function, respectively.

In a renewal risk model, we call $X_i=Y_i-cZ_i$ the net-loss when the $i$th claim comes with distribution $F$
on $(-\infty,\infty)$ for $i\ge1$.
Naturally, the finite-time ruin probability at time $t\ge0$ for some insurance business with initial capital $x\ge0$ is defined by
\begin{eqnarray}\label{104}
\psi(x;t)=P\Big(\max\limits_{0\leq n\leq N(t)}\sum_{i=1}^nX_i>x\Big),
\end{eqnarray}
where $\sum_{i=1}^0X_i$ is defined by $0$.

Further, let random time $\tau$ be a nonnegative random variable satisfying
$$P(\tau>0)>0.$$
Then the random-time ruin probability is defined by
\begin{eqnarray}\label{105}
\psi(x;\tau)=P\Big(\max\limits_{0\leq n\leq N(\tau)}\sum_{i=1}^nX_i>x\Big),
\end{eqnarray}
which is also the main research goal of this paper.

In this paper, we always assume that $\{Y_i,\ i\geq1\}$, $\{Z_i,\ i\geq1\}$ and $\tau$ are mutually independent.

\subsection{\normalsize\bf Main result}


%

First of all, we briefly review some existing results about random-time ruin probability.

For standard renewal risk model, based on the results of finite-time ruin probability,
Corollary 3.2 (1) of Tang \cite{T2004} establishes the following result of random-time ruin probability.
\begin{thm}\label{thm101}
Consider the standard renewal risk model. If $F\in\mathcal{C}$, $EY_1^p<\infty$ for some $p>J_F^++1$ and $E\tau<\infty$,
then it holds that
\begin{eqnarray}\label{106}
\psi(x;\tau)\sim EN(\tau)\overline{G}(x).
\end{eqnarray}
\end{thm}

In order to expand the scope of the distribution class for claim sizes,
Theorem 2.1 of Wang et al. \cite{WGWL2009} provides the following result under the stronger condition (\ref{107}) below.
\begin{thm}\label{thm102}
In the standard renewal risk model, the following assertions are equivalent.

(1) $F\in\mathcal{S}$.

(2) The relation (\ref{106}) holds for any $\tau$ satisfying condition that
\begin{eqnarray}\label{107}
Ee^{r_0N(\tau)}<\infty\ \ \ \ \ \ \ \text{for some}\ \ r_0>0.
\end{eqnarray}

(3) $W\in\mathcal{S}$, where W is the distribution of $\sup_{\{0\le n\le N(\tau)\}}\sum_{i=1}^n(Y_i-cZ_i)$.
\end{thm}

For random-time ruin probabilities in a nonstandard renewal risk model,
building upon existing results for finite-time ruin probabilities,
we reference Theorems 2.1 and 2.2 in Wang et al. \cite{WCWM2012}, or see Theorem \ref{thm103} below.
These results slightly generalize and improve the corresponding results of Tang \cite{T2004}, Leipus and \v{S}iaulys \cite{LS2007},
Kocetova et al. \cite{KLS2009} and Yang et al. \cite{YLSC2011}.

Theorems 2.1 and 2.2 in Wang et al. \cite{WCWM2012} require the following conditions.

\textbf{Condition 1.} The inter-occurrence times $\{Z_i, i\geq1\}$ are NLOD random variables.

\textbf{Condition 2.} The inter-occurrence times $\{Z_i, i\geq1\}$ are WOD random variables
with dominating coefficients $g_{L,H}(n),g_{U,H}(n),n\ge1$.
And there exists a positive and nondecreasing function $g_H(\cdot)$ on $[0,\infty)$ such that $g_H(x)\underline{\uparrow}\infty$,
$x^{-s}g_H(x)\overline{\downarrow}$ for some $0<s<1$, $EZ_1g_H(Z_1)<\infty$ and
$$\max\{g_{L,H}(n),g_{U,H}(n)\}\leq g_H(n)\ \ \ \text{for all}\ \ \ n\geq1.$$
Here, $x^{-s}g_H(x)\overline{\downarrow}$ means that
there exists a finite constant $C>0$ such that
$$x_1^{-s}g_H(x_1)\geq
Cx_2^{-s}g_H(x_2)\quad \mbox{for all}\quad 0\leq x_1<x_2<\infty,$$
and $g_H(x)\underline{\uparrow}\infty$ means that $g^{-1}_H(x)\overline{\downarrow}$.

\textbf{Condition 3.} The inter-occurrence times $\{Z_i, i\geq1\}$ are WOD random variables
with dominating coefficients $g_{L,H}(n),g_{U,H}(n),n\ge1$ such that $EZ_1^p<\infty$
for some $p\ge2$ and when $n\to\infty$,
$$\max\{g_{U,H}(n),g_{L,H}(n)\}=o(n^{b})\ \ \ \ \ \ \ \text{for some}\ b>0.$$

\textbf{Condition 4.} The inter-occurrence times $\{Z_i, i\geq1\}$ are WOD random variables
with dominating coefficients $g_{L,H}(n),g_{U,H}(n),n\ge1$ such that $Ee^{rZ_1}<\infty$ for some $r>0$ and when $n\to\infty$,
$$\max\{g_{U,H}(n),g_{L,H}(n)\}=o(e^{cn})\ \ \ \ \ \ \ \text{for any}\ c>0.$$

\begin{thm}\label{thm103}
In a nonstandard renewal risk model, assume that $Y_i,i\ge1$ are independent identically distributed random variable,
and that $Z_i,i\ge1$ satisfy one of \textbf{Condition 1}-\textbf{Condition 4}.
If $G\in\mathcal{S}^*$ and $E\tau<\infty$, then (\ref{106}) holds.
\end{thm}

In the above three results, the claim sizes are independent of each other.
However, in practice, there is often some dependent relationship between the claim sizes.
For example, in an insurance business, the part of the total claim sizes exceeding a certain limit will be borne by a reinsurance company. In other words, the total claim sizes that the insurance business needs to pay do not exceed this limit.
According to Theorem 2.6 of Jog-Dev and Proschan \cite{JP1983}, see the following proposition,
the relationship between claim sizes is close to NA.
Therefore, it is more reasonable to assume that this relationship is ENOD or WOD.
Similarly, when the operation time $t$ is given, the sum of the claim interval time $\sum_{i=1}^{N(t)}Z_i$ does not exceed $t$.
Therefore, we also have reason to assume that they also have some WOD structure.
\begin{pron}\label{Proposition104}
Let $\xi_i,1\le i\le k$ for some $k\ge1$ be independent and suppose that the condition expectation $Ef(\xi_i,i\in A|\sum_{i\in A}\xi_i)$
is increasing in $\sum_{i\in A}\xi_i$, for every increasing function $f(\cdot)$ and every proper subset A of $\{1,\cdots,k\}$.
Then the condition distribution given $\sum_{i\in A}\xi_i$, is NA almost surely.
\end{pron}

We first give asymptotic estimate of the random time ruin probability for ENUOD claim sizes.
To establish this result, we restrict the distribution class of claim sizes to a narrower scope.
\begin{thm}\label{thm104}
In a nonstandard renewal risk model, let $Y_i,i\ge1$ be ENUOD random variables with dominating coefficient $M_G$,
and let $Z_i,i\ge1$ be WLOD random variables with dominating coefficients $g_{L,H}(n),n\ge1$ satisfying
\begin{eqnarray}\label{109}
g_{L,H}(n)=o(n^{b})\ \ \ \ \ \ \text{for some}\ b>0,\ \ \ \text{as}\ \ \ n\to\infty.
\end{eqnarray}
If $G\in\mathcal{C}$ and
\begin{eqnarray}\label{110}
P(\tau>x)=o\big(\overline{G}(x)\big),
\end{eqnarray}
then $\max\{E\tau,EN(\tau)\}<\infty$ and (\ref{106}) holds.
\end{thm}

Then we expand the dependent structure of the claim sizes from ENOD to WUOD.
To this end, we slightly strengthen the moment condition for $Y_1$.

\begin{thm}\label{thm105}
In a nonstandard renewal risk model, let $Y_i,i\ge1$ be WUOD random variables with dominating coefficients $g_{U,G}(n),n\ge1$,
and let $Z_i,i\ge1$ be WLOD random variables with dominating coefficients $g_{H,G}(n),n\ge1$.
If $G\in\mathcal{C}$, $EY_1^r<\infty$ for some $r>1$, (\ref{109}) and (\ref{110}) are satisfied and
\begin{eqnarray}\label{112}
g_{U,G}(n)=o(n^{r-1})\,\ \ \ \text{as}\ \ \ n\to\infty,
\end{eqnarray}
then $\max\{E\tau,EN(\tau)\}<\infty$ and (\ref{106}) holds.
\end{thm}

\begin{remark}\label{rem101}
In Theorems \ref{thm104} and \ref{thm105}, we need condition (\ref{110}).
According to Proposition \ref{pron102} (3), there is a sufficient condition for (\ref{110}) that,
if $G\in\mathcal{D}$ and $E\tau^p<\infty$ for some $p>J_G^+$, then (\ref{110}) holds.
Therefore, according to Lemma \ref{lem301} below,
(\ref{110}) can be implied by the condition (\ref{107}).
In this way, we can not use the results of finite-time ruin probability like proofs of Theorems \ref{thm101}-\ref{thm103}
to prove Theorems \ref{thm104} and \ref{thm105} in this paper.
We will complete our proof using some results of precise large deviation and renewal theorems
for WOD random variables.
\end{remark}

Let $\tau=t$ for each $t>0$ in the above Theorem \ref{thm104} and Theorem \ref{thm105},
we directly get an asymptotic estimation of the finite-time ruin probability.
However, as mentioned earlier, we hope to give a uniform asymptotic estimation on the operating time $t$ of the finite-time ruin probability, which shows that the initial capital of an insurance business has basically nothing to do with the operation time $t$
in an interval, see, for example, $(0,xg(x)]$ in Theorem \ref{thm106} below.
Specifically, for each pair of $t_2>t_1>0$ and any $0<\varepsilon<1$,
there exists a common threshold $x_0>0$ independent of $t_1$ and $t_2$ such that
$$1-\varepsilon<\psi(x;t_j)EN^{-1}(t_j)\overline{G}^{-1}(x)<1+\varepsilon\ \ \ \ \ \ \text{for all}\ \ \ x\ge x_0,\ \ \ j=1,2.$$

\begin{thm}\label{thm106}
(1) Under the assumptions of Theorem \ref{thm104} but without requiring condition (\ref{110}),
for any positive function $g(\cdot)$ on $[0,\infty)$ satisfying $g(x)\downarrow0$ sufficiently slowly such that
\begin{eqnarray}\label{102111}
x^{-1}g^{-1}(x)\ln x\downarrow0,
\end{eqnarray}
it holds that
\begin{eqnarray}\label{1127}
\lim\sup_{0<t\le xg(x)}\big|\psi(x;t)E^{-1}N(t)\overline{G}^{-1}(x)-1\big|=0.
\end{eqnarray}

(2) Similarly, under the assumptions of Theorem \ref{thm105} without requiring condition (\ref{110}),
the corresponding asymptotic relation (\ref{1127}) also holds.
\end{thm}

\begin{remark}\label{rem102}
(1) The uniform asymptotic estimation for finite-time ruin probability established in Theorem \ref{thm106} constitutes a novel contribution to risk theory, particularly considering the ENOD structure
and even the more generalized WUOD structure imposed on claim sizes.
This result not only shows that the initial capital of an insurance business has essentially nothing
to do with the operation time, but also shows that the interest rate and some related dependence structures in this theorem have no influence on the obtained conclusions under some conditions.

(2) Clearly, if the function $g(\cdot)$ is larger in (\ref{1127}), the range of $t$ is larger.
For example, for any positive function $l(\cdot)$ on $[0,\infty)$ and the function $g(\cdot)$ in Theorem \ref{thm106},
if $l(x)\uparrow\infty$, $l(x)g(x)\downarrow0$, then (\ref{1127}) is also true for function $\widehat{g}(\cdot)=l(\cdot)g(\cdot)$.
In this way, the scope of $t$ is expanded from $0<t\le xg(x)$ to $0<t\le xl(x)g(x)$.
\end{remark}

The paper is organized as follows. We prove Theorem \ref{thm104}-Theorem \ref{thm106} in Section 3, respectively.
To this end, we introduce some results on the precise large deviations
and some renewal theorems for some WOD random variables in Section 2.

\section{\normalsize\bf Some preliminary results}
\setcounter{equation}{0}\setcounter{thm}{0}\setcounter{lemma}{0}\setcounter{remark}{0}\setcounter{exam}{0}\setcounter{pron}{0}
\setcounter{Corol}{0}

In order to prove Theorem \ref{thm104}-Theorem \ref{thm106},
we need some renewal theorems and results of precise large deviation for WOD random variables.

\subsection{\normalsize\bf Two renewal theorems}

In the following two results for the quasi-renewal counting process $N(t)$ generated by WLOD random variables $Z_i,i\ge1$,
the first result comes from Lemma 2.2 of Wang and Cheng \cite{WC2011}.

\begin{thm}\label{thm201}
In a nonstandard renewal risk model,
let $Z_i,\ i\ge1$ be WLOD random variables with dominating coefficients $g_{L,H}(n),\ n\ge1$.
If
\begin{eqnarray}\label{201}
g_{L,H}(n)=o(e^{an})\ \ \ \ \ \ \ \text{for some }\ \ a>0,\ \ \ \text{as}\ \ \ n\to\infty,
\end{eqnarray}
then for each $\delta>0$, there exists $r=r(H)>0$ such that
\begin{eqnarray}\label{202}
Ee^{rN(t)}\textbf{\emph{1}}_{\{N(t)>(1+\delta)\mu_H^{-1}t\}}\to0,\ \ \ \ \ \ \ \text{as}\ \ \ t\to\infty.
\end{eqnarray}
\end{thm}

The second result dues to Lemma 4.2 of Wang et al. \cite{WCWM2012}.

\begin{thm}\label{thm202}
In a nonstandard renewal risk model, let $Z_i,\ i\ge1$ be WLOD random variables with the dominating coefficients $g_{L,H}(n),\ n\ge1$.
If (\ref{109}) is satisfied,
then for each $q\ge1$,
\begin{eqnarray}\label{204}
EN^q(t)\sim\mu_H^{-q}t^p,\ \ \ \ \ \ \ \text{as}\ \ \ t\to\infty.
\end{eqnarray}
\end{thm}

It is clear that condition (\ref{201}) is stronger than condition (\ref{109}).

\subsection{\normalsize\bf Two results on precise large deviations}

The following result is a noncentralized version of Theorem 2.2 of Liu $\cite{L2009}$.

\begin{thm}\label{thm203}
In the above nonstandard renewal risk model, let $Y_i,i\geq 1$ be ENOD random variables with dominating coefficient $M_G$.
If $G\in\mathcal{C}$, then for each $\gamma>\mu_G$,
\begin{eqnarray}\label{205}
\lim\limits_{n\to\infty}\sup\limits_{x\ge\gamma n}\big|P(S^G_n>x)\big(n\overline{G}(x-\mu_G n)\big)^{-1}-1\big|=0.
\end{eqnarray}
\end{thm}
\proof Let $\xi_1=Y_1-\mu_G$ with distribution $V$.
Clearly, $E\xi_1=0$, $P(\xi_1^->x)=o\big(\overline{V}(x)\big)$, $E(\xi_1^-)^r<\infty$ for each $r>1$ and
$$P(S^G_n>x)=P\Big(\sum_{i=1}^n(Y_i-\mu_G)>x-n\mu_G\Big)=P\Big(\sum_{i=1}^n\xi_i>x-n\mu_G\Big).$$
Further, for each $\gamma>\mu_G$, by $x\ge\gamma n$, we have $\gamma_0=\gamma-\mu_G>0$ and
$$x-n\mu_G\ge(\gamma-\mu_G)n=\gamma_0n.$$
In fact, $\gamma_0>0$ is equivalent to $\gamma>\mu_G$.
Therefore, according to Theorem 2.2 of Liu \cite{L2009},
\begin{eqnarray*}
\lim\limits_{n\to\infty}\sup\limits_{x\ge\gamma n}\Big|\frac{P(S^G_n>x)}{n\overline{G}(x-\mu_G n)}-1\Big|
=\lim\limits_{n\to\infty}\sup\limits_{x-n\mu_G\ge\gamma_0 n}\Big|\frac{P(S^G_n-n\mu_G>x-n\mu_G)}{n\overline{G}(x-\mu_G n)}-1\Big|=0,
\end{eqnarray*}
that is (\ref{205})  holds.
$\hspace{\fill}\Box$\\

To prove the second result of precise large deviations, we need the following lemma,
which comes from Proposition 1, Proposition 2 and Theorem 1 of Wang et al. \cite{CWC2016}
and also show that the scope of WOD dependent structure of claim sizes.
Among them, (1) and (2) of this lemma are supplements and adjustments to the above Proposition 2.

\begin{lemma}\label{lem201}
Let $\xi_{i},i\geq1$ be WUOD random variables with a common distribution $V$ on $(-\infty,\infty)$
such that $E\xi_1s_1(\xi_1)<\infty$ for some positive even function $s_1(\cdot)$ on $(-\infty,\infty)$
satisfying $s_1(x)\underline{\uparrow}\infty$.
Then the following three conclusion hold.

(1) There exists a positive even function $s_2(\cdot)$ on $(-\infty,\infty)$ such that
\begin{eqnarray}\label{206}
s_2(x)\underline{\uparrow}\ \infty,\ \ s_2(x)s_1^{-s+1}(x)\overline{\downarrow}\ 0\ \ \ \text{for any}\ \ \ s>1\ \ \
\text{and}\ \ \ E\xi_1s_1(\xi_1)s_2(\xi_1)<\infty.
\end{eqnarray}

(2) Write $s(\cdot)=s_1(\cdot)s_2(\cdot)$ and $V_{s(\cdot)}^+(x)=s^{-1}(x)\int_{0\le y\le x}ys(y)V(dy)+x\overline{V}(x),\ x\ge0$, then
\begin{eqnarray}\label{2c}
V^+_{s(\cdot)}(x)=o\big(s^{-1}_1(x)\big)\ \ \ \ \ \text{and}\ \ \ \ \ s_1^{-l}(x)=o\big(V^+_{s(\cdot)}(x)\big)
\ \ \ \text{for any}\ \ l>1,
\end{eqnarray}

(3) If $E\xi_1=0$ and there is an integer $m\ge1$ such that
$$x^ms^{-1}(x)\underline{\uparrow}\ \infty\ \ \ \ \text{and}\ \ \ x^{m-1}s^{-1}_1(x)\overline{\downarrow}\ 0,$$
then for any $\gamma>0$, $v>0$ and $0<\theta<1$, there exist $x_0=x_0\big(v,\gamma,m,\theta,V,s_1(\cdot),s_2(\cdot)\big)\ge 1$
and $C=C(x_0)>0$ such that for each integer $n\ge1$,
\begin{eqnarray}\label{e1}
P\Big(\sum_{i=1}^{n}\xi_i>x\Big)\le n\overline{V}(vx)+Cs_{U,V}(n)\big(V_{s_{(\cdot)}}^+(vx)\big)^{(1-\theta)v^{-1}}
\ \ \ \ \ \ \ \text{for all}\ \ x\ge \max\{\gamma n,x_0\}.
\end{eqnarray}
\end{lemma}
\proof We only need to prove (1) and (2) of the lemma.

(1) For each integer $n\ge1$, by $E\xi_1s_1(\xi_1)<\infty$, there is a $x_n>0$ large enough such that
$$E\xi_1s_1(\xi_1)\textbf{1}_{\{x_n\le \xi_1<x_{n+1}\}}<n^{-3}.$$
We might as well set up that $n^2<x_n<x_{n+1}$ for all $n\ge1$.
Let $s_3(\cdot)$ be a positive even function on $(-\infty,\infty)$ such that
$$s_3(x)=\textbf{1}_{[0,x_1)}(x)+\sum_{n=1}n\textbf{1}_{[x_n,x_{n+1})}(x)\ \ \ \ \ \ \ \ \ \text{for}\ \ \ x\ge0.$$
Further, for any $s>1$, we take $l_0$ such that $\min\{s,2\}>l_0>1$ and define a positive even function
$s_2(\cdot)$ on $(-\infty,\infty)$ such that
$$s_2(x)=\min\{s_1^{l_0-1}(x),s_3(x)\}\ \ \ \ \ \ \ \ \ \text{for}\ \ \ x\ge0.$$
Then (\ref{206}) holds for the above $s>1$.

(2) Here we only need to prove the second result of (\ref{2c}). By (\ref{206}), we have
$$s_1^l(x)V^+_{s(\cdot)}(x)=s^{l-1}_1(x)s^{-1}_2(x)\int_{0-}^xys(y)V(dy)\to\infty,$$
that is the second result in (\ref{2c}) holds.
\hfill$\Box$

\begin{exam}\label{exam201}
If $E\xi_1^r<\infty$ for some $r\ge1$, then there is an integer $m\ge1$ such that $m\le r<m+1$.
Further, if $r>1$, we first take a positive even function $s_1(\cdot)$ on $(-\infty,\infty)$
such that $s_1(x)=x^{r-1}$ on $(0,\infty)$, then according to Lemma \ref{lem201} (1), there exists a positive even function $s_2(\cdot)$ on $(-\infty,\infty)$ that rises to infinity slowly enough such that (\ref{206}) holds.
And if $r=1$, using the proof of Lemma \ref{lem201} (1), we know that there is a positive even function $s_1(\cdot)$ on $(-\infty,\infty)$ that rises to infinity slowly enough to make $E\xi_1 s_1(\xi_1)<\infty$.
Further, there also exists a positive even function $s_2(\cdot)$ on $(-\infty,\infty)$ such that (\ref{206}) holds.
\hfill$\Box$
\end{exam}

Up to now, precise large deviation results for WUOD random variables,
such as those established in Theorem \ref{thm203} for ENOD random variables, remain unavailable.
In what follows, we present a preliminary estimate.
Despite its limited precision, this result nevertheless enables us to derive an exact asymptotic expression
for the random-time ruin probability.
\begin{thm}\label{thm204}
In a nonstandard renewal risk model, let $Y_i,i\geq 1$ be WUOD random variables
with dominating coefficient $g_{U,G}(n)$ for all $n\ge1$ and $EY_1^r<\infty$ for some $r>1$.

(1) If $G\in\mathcal{C}$ and
\begin{eqnarray}\label{208000}
g_{U,G}(n)=o(n^d)\ \ \ \ \ \text{for some}\ \ \ \ d>0,
\end{eqnarray}
then for each $\gamma>\mu_G$,
\begin{eqnarray}\label{210}
\limsup\limits_{n\to\infty}\sup\limits_{x\ge\gamma n}P(S^G_n>x)\big(n\overline{G}(x)\big)^{-1}
\le\overline{G^*}\big((r-1)(J^+_G+d-1)^{-1}\big)\overline{G_*}^{-1}((1-\mu_G\gamma^{-1})^{-1}).
\end{eqnarray}

(2) If condition (\ref{112}) is satisfied,
then
\begin{eqnarray}\label{2101}
\liminf\limits_{n\to\infty}\inf\limits_{x\ge\gamma n}P(S^G_n>x)\big(n\overline{G}(x)\big)^{-1}\ge1.
\end{eqnarray}
\end{thm}
\proof (1) For any $0<v=v(G)<(r-1)(J^+_G+d-1)^{-1}$,
there is a $0<\theta=\theta(v)<1$ such that
$$(r-1)(1-\theta)v^{-1}-d+1>J^+_G.$$
And take $s_1(x)=x^{r-1},\ x\ge0$. Then according to Proposition \ref{pron102} (3) and the first conclusion in (\ref{2c}),
and by (\ref{208000}), when $n\to\infty$ and $x\ge\gamma n$, we have
\begin{eqnarray}\label{211}
g_{U,G}(n)\big(V_{s(\cdot)}^+(vx)\big)^{\frac{1-\theta}{v}}=o\big(n^ds^{-\frac{1-\theta}{v}}_1(vx)\big)
=o\big(nx^{d-1-\frac{(r-1)(1-\theta)}{v}}\big)=o\big(n\overline{G}(x)\big);
\end{eqnarray}
and by $G\in\mathcal{C}$, we know that
\begin{eqnarray}\label{212}
\overline{G}(vx)\lesssim\overline{G^*}(v)\overline{G}(x).
\end{eqnarray}
Therefore, according to Lemma \ref{lem201}, combining (\ref{e1}), (\ref{211}) and (\ref{212}), for any $\gamma>\mu_G$, we have
\begin{eqnarray*}
&&\limsup\limits_{n\to\infty}\sup\limits_{x\ge\gamma n}P(S^G_n>x)\big(n\overline{G}(x)\big)^{-1}
=\limsup\limits_{n\to\infty}\sup\limits_{x\ge\gamma n}P(S^G_n-n\mu_G>x-n\mu_G)\big(n\overline{G}(x)\big)^{-1}\nonumber\\
&\le&\limsup\limits_{n\to\infty}\sup\limits_{x\ge\gamma n}\overline{G}\big(v(x-n\mu_G)\big)\overline{G}^{-1}(x)\nonumber\\
&\le&\limsup\limits_{n\to\infty}\sup\limits_{x\ge\gamma n}\overline{G^*}(v)\overline{G}(x-n\mu_G)\overline{G}^{-1}(x)\nonumber\\
&\le&\limsup\limits_{n\to\infty}\sup\limits_{x\ge\gamma n}\overline{G^*}(v)\overline{G}\big(x(1-\mu_G\gamma^{-1})\big)\overline{G}^{-1}(x)\nonumber\\
&\le&\overline{G^*}(v)\overline{G_*}^{-1}((1-\mu_G\gamma^{-1})^{-1}).
\end{eqnarray*}
Finally, let $v\uparrow(r-1)(J^+_G+d-1)^{-1}$, (\ref{210}) holds.

(2) From (\ref{112}) and $EY_1^r<\infty$ for some $r>1$, we have
\begin{eqnarray*}
&&\liminf\limits_{n\to\infty}\inf\limits_{x\ge\gamma n}P(S^G_n>x)\big(n\overline{G}(x)\big)^{-1}
\ge\liminf\limits_{n\to\infty}\inf\limits_{x\ge\gamma n}P(\max\{Y_i,1\le i\le n\}>x)\big(n\overline{G}(x)\big)^{-1}\nonumber\\
&\ge&\liminf\limits_{n\to\infty}\inf\limits_{x\ge\gamma n}n\overline{G}(x)\big(1-g_{U,G}(n)n\overline{G}(x)\big)\big(n\overline{G}(x)\big)^{-1}\nonumber\\
&\ge&\limsup\limits_{n\to\infty}\sup\limits_{x\ge\gamma n}\big(1-n^r\overline{G}(x)\big)\nonumber\\
&\ge&\limsup\limits_{n\to\infty}\sup\limits_{x\ge\gamma n}\big(1-\gamma^{-r}x^r\overline{G}(x)\big)\nonumber\\
&=&1,
\end{eqnarray*}
that is (\ref{2101}) holds.
\hfill$\Box$\\

\section{\normalsize\bf Proofs of Theorem \ref{thm104}-Theorem \ref{thm106}}
\setcounter{equation}{0}\setcounter{thm}{0}\setcounter{pron}{0}\setcounter{lemma}{0}\setcounter{Corol}{0}

\subsection{\normalsize\bf Proof of Theorem \ref{thm104}}

Clearly, by $\mu_G<\infty$ and (\ref{110}), we have $E\tau<\infty$.
In order to prove $EN(\tau)<\infty$, we need the following lemma.

\begin{lemma}\label{lem301}
In a nonstandard renewal risk model,
let $Z_i,\ i\ge1$ be WLOD random variables with dominating coefficients $g_{L,H}(n),n\ge1$.
If (\ref{109}) is satisfied, then for each $p\ge1$,
$$E\tau^p<\infty\ \ \ \ \text{if and only if}\ \ \ \ EN^p(\tau)<\infty.$$
\end{lemma}

\proof According to Theorem \ref{thm202}, by (\ref{109}) with some $p\ge1$,
we know that for any $\delta\in(0,1)$, there exists a $t_0=t_0(H,p,\delta)>0$ such that,
for all $t\ge t_0$, it holds that
\begin{eqnarray}\label{301}
(1-\delta)\mu^{-p}_Ht^p<EN^p(t)<(1+\delta)\mu^{-p}_Ht^p.
\end{eqnarray}
If $E\tau^p<\infty$, then by (\ref{301}), we have
\begin{eqnarray*}\label{302}
&&EN^p(\tau)=\Big(\int_{0-}^{t_0}+\int_{t_0}^\infty\Big)EN^p(t)P(\tau\in dt)\nonumber\\
&\le&\int_{0-}^{t_0}EN^p(t_0)P(\tau\in dt)+\int_{t_0}^\infty EN^p(t)P(\tau\in dt)\nonumber\\
&<&(1+\delta)\mu^{-p}_H\int_{0-}^\infty(t_0^p+t^p)P(\tau\in dt)<\infty.
\end{eqnarray*}
Similarly, if $EN^p(\tau)<\infty$ for some $p\ge1$, we can prove $E\tau^p<\infty$.
\hfill$\Box$\\

%

According to Lemma \ref{lem301}, by $E\tau<\infty$, we know that $EN(\tau)<\infty$.

Now, we prove (\ref{106}). To this end, we take any $\gamma>\mu_G$ to split
\begin{eqnarray}\label{311}
&&\psi(x;\tau)=\Big(\int_{0-}^{\gamma^{-1}x}+\int_{\gamma^{-1}x}^\infty\Big)P\big(\max_{1\le k\le N(t)}S^F_k>x\big)P(\tau\in dt)\nonumber\\
&=&I_1(x)+I_2(x),\ \ \ \ \ \ \ \ \ x>0.
\end{eqnarray}

We first deal with $I_2(x)$.  By (\ref{110}) and $G\in\mathcal{C}$, we know that
\begin{eqnarray}\label{313}
I_2(x)\le P(\tau>\gamma^{-1}x)=o\big(\overline{G}(x)\big).
\end{eqnarray}
Therefore, we only need to deal with $I_1(x)$. To this end, for any integer $n_0\ge1$, we again split
\begin{eqnarray}\label{316}
I_1(x)=\int_{0-}^{\gamma^{-1}x}\Big(\sum_{n=1}^{n_0}+\sum_{n=n_0+1}^\infty\Big)
P\big(\max_{1\le k\le n}S^F_k>x,\ N(t)=n\big)P(\tau\in dt)=I_{11}(x)+I_{12}(x).
\end{eqnarray}

To deal with $I_{11}(x)$, we first invoke Lemma 2.3 from Wang et al. \cite{WWG2013}, which was subsequently generalized by Liu et al. \cite{LGW2012} in their Lemma 2.1 to a broader framework extending beyond the WUOD structure.
\begin{lemma}\label{lem303}
Let $\xi_i, i\geq 1$ be nonnegative and WUOD random variables with common distribution $V$ on $[0,\infty)$.
If $V\in\mathscr{L}\cap\mathscr{D}$, then for any fixed $0<a\leq b<\infty$ and each integer $n$, the relation
$$P\left(\sum_{k=1}^n c_k\xi_k>x\right)\sim \sum_{k=1}^nP(c_k\xi_k>x)$$
holds uniformly for all $(c_1,\cdots,c_n)\in[a,b]^n$. 
\end{lemma}
According to Lemma \ref{lem303} , by (\ref{109}) and $EN(\tau)<\infty$,
for the above $n_0$ and any $0<\varepsilon<1$, there exists a $x_0=x_0(G,n_0,\varepsilon)$ large enough, when $x\ge x_0$, we have
\begin{eqnarray}\label{315}
(1-\varepsilon)n\overline{G}(x)\le P(S^G_n>x)\le (1+\varepsilon)n\overline{G}(x)\ \ \ \ \ \text{for all}\ \ \ 1\le n\le n_0
\end{eqnarray}
and
\begin{eqnarray}\label{40000}
\max\{EN(\tau)\textbf{1}_{\{N(\tau)>n_0\}},\ EN(\tau)\textbf{1}_{\{N(\tau)>\gamma^{-1}x\}}\}<2^{-1}\varepsilon EN(\tau).
\end{eqnarray}

For $I_{11}(x)$, on the one hand, by (\ref{315}), for the above $\varepsilon$, $n_0$ and $x_0$, when $x\ge x_0$,
\begin{eqnarray}\label{317}
&&I_{11}(x)\le\int_{0-}^{\gamma^{-1}x}\sum_{n=1}^{n_0}P(S^G_n>x)P\big(N(t)=n\big)P(\tau\in dt)\nonumber\\
&<&(1+\varepsilon)\int_{0-}^{\gamma^{-1}x}\sum_{n=1}^{n_0}nP\big(N(t)=n\big)P(\tau\in dt)\overline{G}(x)\nonumber\\
&\le&(1+\varepsilon)EN(\tau)\overline{G}(x);
\end{eqnarray}
on the other hand, using the proof method of Theorem \ref{thm204} (2),
also according to Lemma \ref{lem303} and by (\ref{40000}), for the above $\varepsilon$ and $n_0$ large enough again,
we take $x_0$ large enough again, then when $x\ge x_0$,
\begin{eqnarray}\label{317-2}
&&I_{11}(x)\geq\int_{0-}^{\gamma^{-1}x}\sum_{n=1}^{n_0}P\big(S_n^F>x,N(t)=n\big)P(\tau\in dt)\nonumber\\
&\ge&\int_{0-}^{\gamma^{-1}x}\sum_{n=1}^{n_0}P\big(S_n^G>x+ct)P\big(N(t)=n\big)P(\tau\in dt)\nonumber\\
&\ge&\int_{0-}^{\gamma^{-1}x}\sum_{n=1}^{n_0}P\big(\max\{Y_i,1\le i\le n\}>x+ct)P\big(N(t)=n\big)P(\tau\in dt)\nonumber\\
&>&(1-\varepsilon)\int_{0-}^{\gamma^{-1}x}\sum_{n=1}^{n_0}n\overline{G}(x+ct)
\big(1-M_G n\overline{G}(x+ct)\big)P\big(N(t)=n\big)P(\tau\in dt)\nonumber\\
&>&(1-\varepsilon)^2\overline{G_*}(1+c\gamma^{-1})EN(\tau)\textbf{1}_{\{N(\tau)\le n_0\}}\textbf{1}_{\{\tau\le\gamma^{-1}x\}}\overline{G}(x)\nonumber\\
&\ge&(1-\varepsilon)^2\overline{G_*}(1+c\gamma^{-1})EN(\tau)(1-\textbf{1}_{\{N(\tau)>n_0\}}
-\textbf{1}_{\{\tau>\gamma^{-1}x\}})\overline{G}(x)\nonumber\\
&\ge&(1-\varepsilon)^3\overline{G_*}(1+c\gamma^{-1})EN(\tau)\overline{G}(x).
\end{eqnarray}
%

Now, we deal with $I_{12}(x)$. To this end, we further split it as follows:
\begin{eqnarray}\label{318}
&&I_{12}(x)=\int_{0-}^\infty\sum_{n=n_0+1}^\infty\big(\textbf{1}_{\{\gamma^{-1}x\ge\max\{t,n\}\}}(t)
+\textbf{1}_{\{0\le t\le{\gamma^{-1}x}<n\}}(t)\big)\nonumber\\
&&\ \ \ \ \ \ \ \ \ \ \ \ \ \ \ \cdot P\big(\max_{1\le k\le n}S_k^F>x,N(t)=n\big)P(\tau\in dt)\nonumber\\
&=&I_{121}(x)+I_{122}(x).
\end{eqnarray}
According to Theorem \ref{thm203}, by $c\mu_H>\mu_G$, $EN(\tau)<\infty$ and (\ref{109}),
for the above $\varepsilon$, we take $n_0$ and $x_0$ large enough again, then
\begin{eqnarray}\label{319}
&&I_{121}(x)\le\int_{0-}^\infty\sum_{n=n_0+1}^\infty P(S_n^G>x)P\big(N(t)=n\big)
\textbf{1}_{\{\gamma^{-1}x\ge\max\{t, n\}\}}(t)P(\tau\in dt)\nonumber\\
&\le&(1+\varepsilon)\int_{0-}^\infty\sum_{n=n_0+1}^\infty n\overline{G}(x-n\mu_G)
P\big(N(t)=n\big)\textbf{1}_{\{\gamma^{-1}x\ge\max\{t, n\}\}}(t)P(\tau\in dt)\nonumber\\
&\le&\frac{(1+\varepsilon)^2\overline{G}(x)}{\overline{G_*}\big((1-\mu_G\gamma^{-1})^{-1}\big)}
\int_{0-}^\infty\sum_{n=n_0+1}^\infty nP\big(N(t)=n\big)\textbf{1}_{\{\gamma^{-1}x\ge\max\{t, n\}\}}(t)P(\tau\in dt)\nonumber\\
&\le&(1+\varepsilon)^2\overline{G_*}^{-1}\big((1-\mu_G\gamma^{-1})^{-1}\big)EN(\tau)\textbf{1}_{\{N(\tau)>n_0\}}\overline{G}(x)\nonumber\\
&<&(1+\varepsilon)^2\varepsilon\overline{G_*}^{-1}\big((1-\mu_G\gamma^{-1})^{-1}\big)EN(\tau)\overline{G}(x).
\end{eqnarray}

To deal with $I_{122}(x)$, we also need a lemma as follows.

\begin{lemma}\label{lem302}
In a nonstandard renewal risk model,
let $Z_i,\ i\ge1$ be WLOD random variables with dominating coefficients $g_{L,H}(n),n\ge1$ satisfying (\ref{109}).
If (\ref{110}) is satisfied, then
\begin{eqnarray}\label{30120}
P\big(N(\tau)>x\big)=o\big(\overline{G}(x)\big).
\end{eqnarray}
\end{lemma}

\proof Clearly, (\ref{201}) can be implied by (\ref{110}).
Thus, according to Theorem \ref{thm201} with any $0<\delta<1$ and some $l=l(H,\delta)>0$, by (\ref{201}),
we know that there exists a $t_0=t_0(H,l,\delta)>0$ such that, it holds that
\begin{eqnarray}\label{3012}
0<Ee^{lN(t)}\textbf{1}_{\{N(t)>(1+\delta)\mu_H^{-1}t\}}<\delta,\ \ \ \ \ \ \ t\ge t_0.
\end{eqnarray}

Now, we take the above $t_0$ and $\mu_Hx(1+\delta)^{-1}$ to split
\begin{eqnarray}\label{3013}
P\big(N(\tau)>x\big)=\Big(\int_{0-}^{t_0}+\int_{t_0}^{\frac{\mu_Hx}{1+\delta}}+\int_{\frac{\mu_Hx}{1+\delta}}^\infty\Big)
P\big(N(t)>x\big)P(\tau\in dt)=\sum_{k=1}^3P_k(x;\tau).
\end{eqnarray}
For $P_1(x;\tau)$, by (\ref{3012}) and $G\in\mathcal{C}$, when $x\ge(1+\delta)\mu_H^{-1}t_0$ for $t_0$ large enough, we have
\begin{eqnarray}\label{30131}
P_1(x;\tau)\le e^{-lx}\int_{0-}^{t_0}Ee^{lN(t_0)}\textbf{1}_{\{N(t_0)>(1+\delta)\mu_H^{-1}t_0\}}
P(\tau\in dt)<\delta\overline{G}(x).
\end{eqnarray}
Similarly,
\begin{eqnarray}\label{30132}
P_2(x;\tau)\le e^{-lx}\int_{t_0}^{\frac{\mu_Hx}{1+\delta}}Ee^{lN(t)}\textbf{1}_{\{N(t)>(1+\delta)\mu_H^{-1}t\}}
P(\tau\in dt)<\delta\overline{G}(x).
\end{eqnarray}
Finally, by (\ref{110}), we take $x_0$ large enough again, then when $x\ge x_0$,
\begin{eqnarray}\label{30133}
P_3(x;\tau)\le P\big(\tau>(1+\delta)^{-1}\mu_Hx\big)<\delta\overline{G}(x).
\end{eqnarray}

Combining (\ref{3013})-(\ref{30133}) and let $\delta\downarrow0$, then (\ref{30120}) holds.
\hfill$\Box$\\

According to Lemma \ref{lem302}, by (\ref{109}) and $EN(\tau)<\infty$,
for the above $\varepsilon$, we take $x_0$ large enough again, then when $x\ge x_0$,
\begin{eqnarray}\label{435}
&&I_{122}(x)\le\int_{0-}^\infty\sum_{n=n_0+1}^\infty P\big(N(t)=n\big)\textbf{1}_{\{0\le t\le\gamma^{-1}x<n\}}(t)P(\tau\in dt)\nonumber\\
&=&\int_{0-}^\infty\sum_{n=n_0+1}^\infty E\textbf{1}_{\{N(t)=n,0\le t\le\gamma^{-1}x<n\}}(t)P(\tau\in dt)\nonumber\\
&=&\int_{0-}^\infty\sum_{n=n_0+1}^\infty E\textbf{1}_{\{N(t)=n,0\le t\le\gamma^{-1}x<N(t)\}}(t)P(\tau\in dt)\nonumber\\
&=&\int_{0-}^\infty P\big(N(t)>n_0,0\le t\le\gamma^{-1}x<N(t)\big)P(\tau\in dt)\nonumber\\
&\le&P\big(N(\tau)>\gamma^{-1}x\big)\nonumber\\
&<&\varepsilon EN(\tau)\overline{G}(x).
\end{eqnarray}

Combining (\ref{311})-(\ref{316}), (\ref{317})-(\ref{319}) and (\ref{435}), writing
$$b_1(\gamma)=\overline{G_*}(1+c\gamma^{-1})\ \ \ \text{and}\ \ \
b_2(\gamma)=\overline{G}^{-1}_*\Big(\big(1+(c-\mu_H^{-1}\mu_G)\gamma^{-1}\big)^{-1}\Big),$$
then we have
\begin{eqnarray}\label{4120001}
(1-\varepsilon)^3b_1(\gamma)\le\liminf\frac{\psi(x;\tau)}{EN(\tau)\overline{G}(x)}
\le\limsup\frac{\psi(x;\tau)}{EN(\tau)\overline{G}(x)}\le(1+\varepsilon)^3b_2(\gamma).
\end{eqnarray}
Let $\gamma\uparrow\infty$ and $\varepsilon\downarrow0$ in (\ref{4120001}), then by $G\in\mathcal{C}$, we know that (\ref{106}) holds.
$\hspace{\fill}\Box$

\subsection{\normalsize\bf Proof of Theorem \ref{thm105}}

We prove this theorem following the same approach as in Theorem \ref{thm104}, adopting the corresponding notation throughout the proof.
In this way, we only need to deal with the upper bound of $I_{121}(x;t)$.

According to Theorem \ref{thm204}, by $EY_1^r<\infty$ for some $r>1$, (\ref{109}), $EN(\tau)<\infty$ and $c\mu_H>\mu_G$, for the any $0<\varepsilon<1$ and $n_0$ large enough, we have
\begin{eqnarray}\label{3190}
&&I_{121}(x)\le\int_{0-}^\infty\sum_{n=n_0+1}^\infty P(S_n^G>x)P\big(N(t)=n\big)
\textbf{1}_{\{\gamma^{-1}x\ge\max\{t, n\}\}}(t)P(\tau\in dt)\nonumber\\
&\le&\frac{(1+\varepsilon)\overline{G^*}\big(\frac{r-1}{J^+_G+d-1}\big)\overline{G}(x)}{\overline{G_*}\big((1-\mu_G\gamma^{-1})^{-1}\big)}
\int_{0-}^\infty\sum_{n=n_0+1}^\infty nP\big(N(t)=n\big)\textbf{1}_{\{\gamma^{-1}x\ge\max\{t, n\}\}}(t)P(\tau\in dt)\nonumber\\
&\le&(1+\varepsilon)\overline{G^*}\big((r-1)(J^+_G+d-1)^{-1}\big)\overline{G_*}^{-1}\big((1-\mu_G\gamma^{-1})^{-1}\big)
EN(\tau)\textbf{1}_{\{N(\tau)>n_0\}}\overline{G}(x)\nonumber\\
&<&\varepsilon(1+\varepsilon)\overline{G^*}\big((r-1)(J^+_G+d-1)^{-1}\big)
\overline{G_*}^{-1}\big((1-\mu_G\gamma^{-1})^{-1}\big)EN(\tau)\overline{G}(x).
\end{eqnarray}

Combining the proof of Theorem \ref{thm104} and (\ref{3190}), we know that (\ref{106}) holds.
$\hspace{\fill}\Box$

\subsection{\normalsize\bf Proof of Theorem \ref{thm106}}

(1) One the one hand, we prove that
\begin{eqnarray}\label{3201}
\liminf\inf_{0<t\le xg_1(x)}\psi(x;t)E^{-1}N(t)\overline{G}^{-1}(x)\ge1
\end{eqnarray}
for any positive function $g_1(\cdot)$ on $[0,\infty)$ that $g_1(x)\to0$ sufficiently slowly such that at least $xg_1(x)\to\infty$.

To this end, we need the following lemma, which requires no assumptions on the distribution $G$ of the claim size .

\begin{lemma}\label{lem304}
(1) Under the assumptions of Theorem \ref{thm104} but without requiring conditions (\ref{110}) and $G\in\mathcal{C}$, we have
\begin{eqnarray}\label{32010}
\lim\inf_{t>0}\psi(x;t)E^{-1}N(t)\overline{G}^{-1}(x+ct)\ge1.
\end{eqnarray}

(2) Under the assumptions of Theorem \ref{thm105} but without requiring conditions (\ref{110}) and $G\in\mathcal{C}$,
(\ref{32010}) still holds.
\end{lemma}

\proof (1) For any $0<\varepsilon,\delta<1$, by $\mu_G<\infty$, there exists a $x_{0,1}=x_{0,1}(G,H,\varepsilon,\delta)>0$ large enough such that
$$M_G (1+\delta)\mu_H^{-1}x\overline{G}(x)<\varepsilon\ \ \ \ \ \ \text{for all}\ \ x\ge x_{0,1}.$$
In addition, for all $t>0$, let $s(\cdot;t)$ be a positive function on $[0,\infty)$ defined as
$$s(x;t)=(1+\delta)\mu_H^{-1}\max\{x,t\}.$$
Thus, when $x\ge x_{0,1}$, for all $t>0$, it holds that
\begin{eqnarray}\label{32000}
&&\psi(x;t)\ge\sum_{1\le n\le s(x;t)}P\big(S_n^G>x+ct)P\big(N(t)=n\big)\nonumber\\
&\ge&\sum_{1\le n\le s(x;t)}P\big(\max\{Y_i,1\le i\le n\}>x+ct)P\big(N(t)=n\big)\nonumber\\
&\ge&\sum_{1\le n\le s(x;t)}n\overline{G}(x+ct)\big(1-M_G n\overline{G}(x+ct)\big)P\big(N(t)=n\big)\nonumber\\
&\ge&\sum_{1\le n\le s(x;t)}n\overline{G}(x+ct)
\big(1-M_G s(x;t)\overline{G}(x+ct)\big)P\big(N(t)=n\big).
\end{eqnarray}

Now, we find the asymptotic lower bound of $\psi(x;t)$ in three cases: $t>x$,
$t_0<t\le x$ and $0<t\le t_0$ for any fixed $t_0>0$ satisfying
\begin{eqnarray}\label{32200}
0<P(Z_1\le t_0)<1.
\end{eqnarray}

In the case that $t_0<t\le x$, for the above $0<\delta,\varepsilon<1$, according to Theorem \ref{thm201} and by (\ref{109}),
we can take and $x_{0,2}\ge x_{0,1}$ large enough such that
$$EN(x)\textbf{1}_{\{N(x)>s(x;t)\}}E^{-1}N(t_0)=EN(x)\textbf{1}_{\{N(x)>(1+\delta)\mu_H^{-1}x\}}E^{-1}N(t_0)<\varepsilon
\ \ \ \ \text{for all}\ \ x\ge x_{0,2}.$$
Thus, for the above $\varepsilon,\delta$ and $x_{0,2}$, when $x\ge x_{0,2}$, by (\ref{32000}),
$\mu_H<\infty$ and $x\ge t$, it holds that
\begin{eqnarray}\label{32001}
&&\psi(x;t)>(1-\varepsilon)\overline{G}(x+ct)EN(t)\big(1-EN(t)\textbf{1}_{\{N(t)>(1+\delta)\mu_H^{-1}x\}}E^{-1}N(t)\big)\nonumber\\
&\ge&(1-\varepsilon)\overline{G}(x+ct)EN(t)\big(1-EN(x)\textbf{1}_{\{N(x)>(1+\delta)\mu_H^{-1}x\}}E^{-1}N(t_0)\big)\nonumber\\
&>&(1-\varepsilon)^2EN(t)\overline{G}(x+ct).
\end{eqnarray}

In the case that $0<t\le t_0$, for the above $0<\varepsilon,\delta<1$, by (\ref{32200}),
we also take $x_{0,3}\ge x_{0,2}$ large enough such that
$$\sum_{n>s(x;t)}n^{b+1}P^{n-1}(Z_1\le t_0)=\sum_{n>(1+\delta)\mu_H^{-1}x}n^{b+1}P^{n-1}(Z_1\le t_0)<\varepsilon
\ \ \ \ \ \ \ \ \text{for all}\ \ x\ge x_{0,3}.$$
Thus, by (\ref{109}), (\ref{32000}), $\mu_H<\infty$ and the fact that
\begin{eqnarray}\label{32100}
P(Z_1\le t)\le EN(t)\ \ \ \ \ \ \text{for all}\ \ t>0,
\end{eqnarray}
when $x\ge x_{0,3}$, it holds that
\begin{eqnarray}\label{3210000}
&&\psi(x;t)>(1-\varepsilon)\overline{G}(x+ct)\sum_{1\le n\le (1+\delta)\mu_H^{-1}x}nP\big(N(t)=n\big)\nonumber\\
&\ge&(1-\varepsilon)\overline{G}(x+ct)\Big(EN(t)-\sum_{n>(1+\delta)\mu_H^{-1}x}nP\big(N(t)=n\big)\Big)\nonumber\\
&\ge&(1-\varepsilon)\overline{G}(x+ct)\Big(EN(t)-\sum_{n>(1+\delta)\mu_H^{-1}x}nP(Z_1\le t,\cdots,Z_n\le t)\Big)\nonumber\\
&\ge&(1-\varepsilon)\overline{G}(x+ct)\Big(EN(t)-\sum_{n>(1+\delta)\mu_H^{-1}x}ng_{L,H}(n)P^{n-1}(Z_1\le t)EN(t)\Big)\nonumber\\
&\ge&(1-\varepsilon)\overline{G}(x+ct)EN(t)\Big(1-\sum_{n>(1+\delta)\mu_H^{-1}x}n^{b+1}P^{n-1}(Z_1\le t_0)\Big)\nonumber\\
&\ge&(1-\varepsilon)^2EN(t)\overline{G}(x+ct).
\end{eqnarray}

In the case that $t>x$, for the above $0<\delta,\varepsilon<1$, according to Theorem \ref{thm201} and by (\ref{109}),
we can take $x_{0,4}\ge x_{0,3}$ large enough such that $EN(t)\ge1$ and
$$EN(t)\textbf{1}_{\{N(t)>s(x;t)\}}E^{-1}N(t)=EN(t)\textbf{1}_{\{N(t)>(1+\delta)\mu_H^{-1}t\}}E^{-1}N(t)<\varepsilon
\ \ \ \ \text{for all}\ \ t>x\ge x_{0,4}.$$
Thus, for the above $\varepsilon,\delta$ and $x_{0,4}$, when $x\ge x_{0,4}$, by (\ref{32000}) and $t>x$, it holds that
\begin{eqnarray}\label{320012}
&&\psi(x;t)\ge\sum_{1\le n\le(1+\delta)\mu_H^{-1}t}n\overline{G}(x+ct)
\big(1-M_G(1+\delta)\mu_H^{-1}t\overline{G}(x+ct)\big)P\big(N(t)=n\big)\nonumber\\
&\ge&(1-\varepsilon)\overline{G}(x+ct)EN(t)\big(1-EN(t)\textbf{1}_{\{N(t)>(1+\delta)\mu_H^{-1}t\}}E^{-1}N(t)\big)\nonumber\\
&\ge&(1-\varepsilon)^2EN(t)\overline{G}(x+ct).
\end{eqnarray}

Combining (\ref{32000}), (\ref{3210000}) and (\ref{320012}), by the arbitrariness of $\varepsilon$, we know that (\ref{32010}) holds.
\\

(2) We omit the similar proof.
$\hspace{\fill}\Box$\\

Now, we prove (\ref{3201}). By $G\in\mathcal{C}$ and $g_1(x)\to0$, it holds that
\begin{eqnarray}\label{32011}
\lim\inf_{0<t\le xg_1(x)}\overline{G}(x+ct)\overline{G}^{-1}(x)\ge\lim\overline{G}\big(x(1+cg_1(x))\big)\overline{G}^{-1}(x)=1.
\end{eqnarray}
Thus, by the proof of Lemma \ref{lem304}, (\ref{3201}) holds.
\\

On the other hand, we prove that
\begin{eqnarray}\label{3202}
\limsup\sup_{0<t\le xg_2(x)}\psi(x;t)E^{-1}N(t)\overline{G}^{-1}(x)\le1
\end{eqnarray}
for any positive function $g_2(\cdot)$ on $[0,\infty)$ satisfying $g_2(x)\downarrow0$ sufficiently slowly
such that (\ref{102111}) is satisfied.

To this end, we split the sum $\psi(x;t)$ using an integer $n_0$ (to be determined later) as follows:
\begin{eqnarray}\label{321}
\psi(x;t)=\Big(\sum_{n=1}^{n_0}+\sum_{n=n_0+1}^\infty\Big)P\big(\max_{1\le k\le n}S^F_k>x,\ N(t)=n\big)=\psi_{1}(x;t)+\psi_{2}(x;t).
\end{eqnarray}

For $\psi_{1}(x;t)$ and any $0<\varepsilon<1$, when $t>0$, according to Lemma \ref{lem303}, we have
\begin{eqnarray}\label{322}
&&\psi_{1}(x;t)\le\sum_{n=1}^{n_0}P(S^G_n>x)P\big(N(t)=n\big)\nonumber\\
&<&(1+\varepsilon)\sum_{n=1}^{n_0}n\overline{G}(x)P\big(N(t)=n\big)\nonumber\\
&=&(1+\varepsilon)EN(t)\textbf{1}_{\{N(t)\le n_0\}}\overline{G}(x).
\end{eqnarray}

For $\psi_{2}(x;t)$, by any $\delta>0$ and the positive function $g_2(\cdot)$ on $[0,\infty)$ satisfying (\ref{102111}),
we continue to split
\begin{eqnarray}\label{323}
&&\psi_{2}(x;t)=\Big(\sum_{n_0+1\le n\le(1+\delta)\mu_H^{-1}xg_2(x)}+\sum_{n>(1+\delta)\mu_H^{-1}xg_2(x)}\Big)
P\big(\max_{1\le k\le n}S_k^F>x,N(t)=n\big)\nonumber\\
&=&\psi_{21}(x;t)+\psi_{22}(x;t),\ \ \ \ \ \ \ \ \ \ \ \ \ \ \ \ \ x\ge0.
\end{eqnarray}

About $\psi_{21}(x;t)$, for the above $\ \delta$ and any $\gamma>\mu_G$,
by $g_2(x)\to0$, we take $x_{0,1}=x_{0,1}(H,g_2(\cdot),\delta)$ large enough such that
$$(1+\delta)\mu_H^{-1}g_2^{-1}(x)>\gamma,\ \ \ \ \ \ x\ge x_{0,1}.$$
Thus, when $x\ge x_{0,1}$ and $n_0+1\le n\le (1+\delta)\mu_H^{-1}xg_2(x)$,
$$x\ge(1+\delta)\mu_H^{-1}g_2^{-1}(x)n>\gamma n.$$
Further, according to Theorem \ref{thm203}, by (\ref{112}) for $r>1$ and $G\in\mathcal{C}$, for the above $\varepsilon$,
we choose $n_0=n_0(G,\varepsilon)$ and $x_{0,2}=x_{0,2}(n_0)\ge x_{0,1}$ sufficiently large such that,
when $n\ge n_0$ and $x\ge x_{0,2}$,
$$P(S_n^G>x)<(1+\varepsilon)\overline{G}(x-n\mu_G)$$
and
$$\overline{G}\big(x(1-\mu_G(1+\delta)\mu_H^{-1}g_2(x))\big)<(1+\varepsilon)\overline{G}(x).$$
Therefore, for each $t>0$, we know that
\begin{eqnarray}\label{324}
&&\psi_{21}(x;t)\le\sum_{n_0+1\le n\le (1+\delta)\mu_H^{-1}xg_2(x)}P(S_n^G>x)P\big(N(t)=n\big)\nonumber\\
&\le&(1+\varepsilon)\sum_{n_0+1\le n\le (1+\delta)\mu_H^{-1}xg_2(x)} n\overline{G}(x-n\mu_G)P\big(N(t)=n\big)\nonumber\\
&\le&(1+\varepsilon)\sum_{n_0+1\le n\le (1+\delta)\mu_H^{-1}xg_2(x)} n\overline{G}\big(x(1-\mu_G(1+\delta)\mu_H^{-1}g_2(x))\big)P\big(N(t)=n\big)\nonumber\\
&\le&(1+\varepsilon)^2\overline{G}(x)\sum_{n_0+1\le n\le (1+\delta)\mu_H^{-1}xg_2(x)}nP\big(N(t)=n\big)\nonumber\\
&\le&(1+\varepsilon)^2EN(t)\textbf{1}_{\{N(t)>n_0\}}\overline{G}(x).
\end{eqnarray}

We deal with $\psi_{22}(x;t)$ in the following two cases that
$$0<t\le f(x)\ \ \ \text{and}\ \ \ f(x)<t\le xg_2(x),$$
where $f(\cdot)$ is a positive function on $[0,\infty)$ defined as
$$f(x)=\textbf{1}_{[0,\widehat{x})}(x)+\max\{y:\overline{H}(y)\ge\mu_H(p+b+1)x^{-1}g_2^{-1}(x)\ln x\}\textbf{1}_{[\widehat{x},\infty)}(x)$$
for $b$ in (\ref{109}), some $p>J_G^+$ and
$$\widehat{x}=\widehat{x}(H,g_2(\cdot),p,b)=\min\{x:\mu_H(p+b+1)x^{-1}g_2^{-1}(x)\ln x\le1\}.$$

Here, we make some explanations for function $f(\cdot)$.
Clearly, by (\ref{102111}), we have
$$\mu_H(p+b+1)g_2^{-1}(x)x^{-1}\ln x=:g_0^{-1}(x)x^{-1}\ln x\downarrow0.$$
Consequently, both $f(x)$ and $\widehat{x}$ are well-defined satisfying
\begin{eqnarray}\label{3203}
f(x)\uparrow\infty\ \ \ \ \text{and}\ \ \ \ \overline{H}\big(f(x)\big)\ge g_0^{-1}(x)x^{-1}\ln x,\ \ \ \ \ \ \ x\ge \widehat{x}.
\end{eqnarray}
We are also happy to give some concrete examples of the function $f(\cdot)$.
\begin{exam}\label{exam101}
For some $\alpha>1$, if
$$\overline{H}(x)=\emph{\textbf{1}}_{(-\infty,1)}(x)+x^{-\alpha}\emph{\textbf{1}}_{[1,\infty)}(x),
\ \ \ \ \ \ \ \ x\in(-\infty,\infty),$$
then we can take $f(\cdot)$ on $[0,\infty)$ such that
$$f(x)=\emph{\textbf{1}}_{[0,e)}(x)+x^{\frac{1}{\alpha}}(\ln x)^{-\frac{1}{\alpha}} g_0^{\frac{1}{\alpha}}(x)\emph{\textbf{1}}_{[e,\infty)}(x),\ \ \ \ \ \ \ \ x\ge0.$$

If
$$\overline{H}(x)=\emph{\textbf{1}}_{(-\infty,0)}(x)+e^{-x^{1/2}}\emph{\textbf{1}}_{[0,\infty)}(x),
\ \ \ \ \ \ \ \ x\in(-\infty,\infty),$$
then we can take $f(\cdot)$ on $[0,\infty)$ such that
$$f(x)=\emph{\textbf{1}}_{[0,e^e)}(x)+\big(\ln x-\ln\ln x+\ln g_0(x)\big)^2\emph{\textbf{1}}_{[e^e,\infty)}(x);$$

If
$$\overline{H}(x)=\emph{\textbf{1}}_{(-\infty,0)}(x)+e^{-x}\emph{\textbf{1}}_{[0,\infty)}(x),\ \ \ \ \ \ \ \ x\in(-\infty,\infty),$$
then we take $f(\cdot)$ on $[0,\infty)$ such that
$$f(x)=\emph{\textbf{1}}_{[0,e^e)}(x)+\big(\ln x-\ln\ln x+\ln g_0(x)\big)\emph{\textbf{1}}_{[e^e,\infty)}(x).$$

For any positive and increasing function $l(\cdot)$ on $[0,\infty)$, the inequality $f(xl(x))\ge f(x)$ holds trivially.
Specially, in the first example, if $l(x)\uparrow\infty$, then
$$f(x)=o\big(f(xl(x))\big).$$
\end{exam}

The following proposition shows that it is reasonable to assume that $f(x)<t\le xg_2(x)$.
\begin{pron}\label{pron105}
For the above functions $f(\cdot)$ and $g_2(\cdot)$, it holds that
\begin{eqnarray}\label{116}
f(x)=o\big(xg_2(x)\big).
\end{eqnarray}
\end{pron}

\proof Because $EZ_1<\infty$, thus $x\overline{H}(x)\to 0$. Hence, by (\ref{3203}),
we can get that
$$f(x)\mu_H(p+b+1)x^{-1}g_2^{-1}(x)\ln x\le f(x)\overline{H}(f(x))\to 0,$$
which implies that (\ref{116}) holds.
$\hspace{\fill}\Box$\\

For the above $0<\varepsilon,\delta<1$, by (\ref{3203}),
we can take $x_{0,3}\ge x_{0,2}$ large enough such that, when $x\ge x_{0,3}$,
\begin{eqnarray*}
&&\int_{(1+\delta)\mu_H^{-1}xg_2(x)}^\infty y^be^{-\overline{H}(f(x))y}dy
=-\overline{H}^{-1}(f(x))\int_{(1+\delta)\mu_H^{-1}xg_2(x)}^\infty y^bde^{-\overline{H}(f(x))y}\nonumber\\
&=&-\overline{H}^{-1}(f(x))\Big(y^be^{-\overline{H}(f(x))y}\big|_{(1+\delta)\mu_H^{-1}xg_2(x)}^\infty
-b\int_{(1+\delta)\mu_H^{-1}xg_2(x)}^\infty y^{b-1}e^{-\overline{H}(f(x))y}dy\Big)\nonumber\\
&\le&\overline{H}^{-1}(f(x))\Big((1+\delta)^b\mu_H^{-b}x^bg_2^b(x)e^{-(1+\delta)\mu_H^{-1}\overline{H}(f(x))xg_2(x)}\nonumber\\
&&\ \ \ \ \ \ \ \ \ \ +b(1+\delta)^{-1}\mu_Hx^{-1}g^{-1}(x)\int_{(1+\delta)\mu_H^{-1}xg_2(x)}^\infty y^{b}e^{-\overline{H}(f(x))y}dy\Big)\nonumber\\
&\le&\frac{(1+\delta)^bx^{b+1}g_2^{b+1}(x)}{\mu_H^{b+1}(p+b+1)e^{(1+\delta)\overline{H}(f(x))\mu_H^{-1}lxg_2(x)}\ln x}
+\frac{bl^{-1}}{(1+\delta)(p+b+1)\ln lx}\int_{\frac{(1+\delta)lxg_2(x)}{\mu_H}}^\infty\frac{y^{b}}{e^{\overline{H}(f(x))y}}dy\nonumber\\
&\le&\frac{(1+\delta)^bx^{b+1}g_2^{b+1}(x)}{\mu_H^{b+1}(p+b+1)e^{\overline{H}(f(x))\mu_H^{-1}xg_2(x)}\ln x}
+\varepsilon\int_{(1+\delta)\mu_H^{-1}xg_2(x)}^\infty\frac{y^{b}}{e^{\overline{H}(f(x))y}}dy,
\end{eqnarray*}
which implies that
\begin{eqnarray}\label{32800}
\int_{(1+\delta)\mu_H^{-1}xg_2(x)}^\infty y^be^{-\overline{H}(f(x))y}dy
\le\frac{(1+\delta)^bx^{b+1}g_2^{b+1}(x)}
{(p+b+1)(1-\varepsilon)\mu_H^{b+1}e^{\overline{H}(f(x))\mu_H^{-1}xg_2(x)}}.
\end{eqnarray}
Thus, according to Theorem \ref{thm203} and Proposition \ref{pron102} (3) again,
by (\ref{32100}), (\ref{109}), (\ref{3203}) and (\ref{32800}),
for the above $0<\delta,\varepsilon<1$, we take $x_{04}\ge x_{03}$ large enough such that,
when $0<t\le f(x)$,
\begin{eqnarray}\label{325}
&&\psi_{22}(x;t)\le\sum_{n>(1+\delta)\mu_H^{-1}xg_2(x)}P\big(N(t)=n\big)\nonumber\\
&\le&\sum_{n>(1+\delta)\mu_H^{-1}xg_2(x)}P(Z_1\le t,\cdots,Z_n\le t)\nonumber\\
&\le&\sum_{n>(1+\delta)\mu_H^{-1}xg_2(x)}g_{L,H}(n)P^n(Z_1\le t)E^{-1}N(t)EN(t)\nonumber\\
&\le&\sum_{n>(1+\delta)\mu_H^{-1}xg_2(x)}n^b\big(1-\overline{H}(t)\big)^{n-1}EN(t)\nonumber\\
&\le&\sum_{n>(1+\delta)\mu_H^{-1}xg_2(x)}\int_n^{n+1}y^b\big(1-\overline{H}(t)\big)^{y-2}dyEN(t)\nonumber\\
&\le&\int_{(1+\delta)\mu_H^{-1}xg_2(x)}^\infty y^b\big(1-\overline{H}(f(x))\big)^{\overline{H}^{-1}(f(x))\overline{H}(f(x))(y-2)}dyEN(t)\nonumber\\
&\le&2\int_{(1+\delta)\mu_H^{-1}xg_2(x)}^\infty y^be^{-\overline{H}(f(x))y}dyEN(t)\nonumber\\
&\le&2(1-\varepsilon)^{-1}(1+\delta)^{b}\mu_H^{-b-1}
(p+b+1)^{-1}x^{b+1}g_2^{b+1}(x)e^{-\overline{H}(f(x))\mu_H^{-1}xg_2(x)}EN(t)\nonumber\\
&\le&x^{b+1}e^{-(p+b+1)\ln x}EN(t)\nonumber\\
&=&x^{-p}EN(t)\nonumber\\
&<&\varepsilon EN(t)\overline{G}(x)\ \ \ \ \ \ \ \ \ \ \text{for all}\ \ x\ge x_{04}.
\end{eqnarray}

In the case that $f(x)<t\le xg_2(x)$, according to Theorem \ref{thm201} and by (\ref{109}) which implies (\ref{201}),
for the above $0<\delta,\varepsilon<1$ and any $p>J_H^+$,
there exists a $x_{0,5}\ge x_{0,4}$ sufficiently large such that, when $x\ge x_{0,5}$, $EN(t)\ge EN(f(x))\ge1$,
$$g_2(x)\ge r^{-1}(1+\delta)^{-1}\mu_Hpx^{-1}\ln x\ \ \ \text{and}\ \ \ Ee^{rN(xg_2(x))}\textbf{1}_{\{N(xg_2(x))>(1+\delta)\mu_H^{-1}xg_2(x)\}}<\varepsilon.$$
Therefore, for the above $0<\delta,\varepsilon<1$, we have
\begin{eqnarray}\label{333}
&&\psi_{22}(x;t)\le P\big(N(xg_2(x))>(1+\delta)\mu_H^{-1}xg_2(x)\big)\nonumber\\
&\le&e^{-r(1+\delta)\mu_H^{-1}xg_2(x)}Ee^{rN(xg_2(x))}\textbf{1}_{\{N(xg_2(x))>(1+\delta)\mu_H^{-1}xg_2(x)\}}\nonumber\\
&\le&\varepsilon e^{-p\ln x}\nonumber\\
&\le&\varepsilon x^{-p}EN(t)\nonumber\\
&\le&\varepsilon EN(t)\overline{G}(x).
\end{eqnarray}

From (\ref{321})-(\ref{324}), (\ref{325}), (\ref{333}) and the arbitrariness of $\varepsilon$,
we know that (\ref{3202}) holds.

Finally, let $g(\cdot)=g_2(\cdot)$, then $g(\cdot)$ also satisfies the conditions of Theorem \ref{thm106} (1).
Therefore, (\ref{1127}) holds by (\ref{3201}) and (\ref{3202}).
\\

(2) We omit the similar proof.
$\hspace{\fill}\Box$\\

\end{document}